\documentclass[letterpaper, 10 pt, conference]{ieeeconf}
\IEEEoverridecommandlockouts
\usepackage{graphicx}
\usepackage{caption}
\usepackage{subcaption}
\usepackage{amsmath}
\usepackage{mathtools}
\usepackage{dsfont}
\usepackage[utf8]{inputenc}
\usepackage{hyperref}
\usepackage{flushend}

\newtheorem{theorem}{Theorem}
\newtheorem{definition}{Definition}
\newtheorem{assumption}{Assumption}
\newtheorem{lemma}{Lemma}
\newtheorem{corollary}{Corollary}

\title{\LARGE \bf
Guaranteed Cost Approach for Robust Model Predictive Control of Uncertain Linear Systems
}

\author{Carlos M. Massera$^{1}$, Marco H. Terra$^{2}$ and Denis F. Wolf$^{1}$
\thanks{This work was supported by FAPESP grant 2013/24542-7.}
\thanks{$^{1}$Institute of Mathmatics and Computer Science, University of São Paulo, Avenida Trabalhador São-carlense, 400, São Carlos, Brazil {\tt\small massera,denis@icmc.usp.br}}%
\thanks{$^{2}$São Carlos School of Engineering, University of São Paulo, Avenida Trabalhador São-carlense, 400, São Carlos, Brazil {\tt\small terra@sc.usp.br}}%
}

\begin{document}
\maketitle
\thispagestyle{empty}
\pagestyle{empty}

\begin{abstract}
In this paper we propose a constrained guaranteed cost robust model predictive controller (GCMPC) for uncertain discrete time systems. 
This controller was developed based on a quadratic cost functional and guarantee  robustness with respect to quadratically bound uncertainties. Such a class of problems is currently intractable by Min-Max Robust Model Predictive Controllers without polytopic approximations of the uncertainties. 
The proposed technique is computationally more efficient then an enumeration-based approach and requires only a Quadratically Constrained Quadratic Problem (QCQP) optimization, whereas LMI-based GCMPC approaches require a Semi-Definite Programming (SDP) optimization. 
\end{abstract}

\section{INTRODUCTION}

Model Predictive Control (MPC) is a class of optimization-based control algorithms that uses an explicit model of the controlled system to predict its future states \cite{badgwell2015model}. 
The technique has been used in several different areas, such as refineries, food processing plants, mining, aerospace and automotive control \cite{qin2003survey}. MPC objective is to minimize a cost functional and maintain the system states and control inputs within a feasible set. The system dynamics are usually assumed to be known. Therefore, model mismatches or external disturbances are not considered. However, the disregard of such uncertainties may lead to poor closed-loop performance and the violation of state and control input constraints \cite{rawlings1994nonlinear}.

Robust Model Predictive Control (RMPC) addresses the poor closed-loop performance subject to uncertainties \cite{bemporad1999robust}. The main RMPC approach is based on Min-Max optimization. Its objective is to obtain a control input sequence that minimizes a cost functional and to guarantee the feasibility of system states and control inputs when the system is subject to the worst-case disturbance. Most Min-Max Model Predictive Controllers (MMMPC) are based on different assumptions and approximations, therefore comparisons among different approaches are difficult. However, most of them are classified into two categories \cite{kerrigan2004feedback}, namely open-loop and feedback MMMPC. The open-loop MMMPC is similar to the nominal MPC, in which a control input sequence is optimized with respect to the open-loop system. It is considered conservative for actual applications. On the other hand, feedback MMMPC optimizes a control input sequence with respect to a disturbance-rejecting closed-loop system, which addresses the conservativeness of the open-loop approach.

MMMPC schemes are defined by exact or approximate solutions to $l1$-norm, $l2$-norm or $l\infty$-norm cost functionals subject to polytopic bounded uncertainties (e.g. bound by $l1$-norm or $l\infty$-norm) \cite{lofberg2003approximations, scokaert1998min, genceli1993robust}. Both $l1$-norm and $l\infty$-norm MMMPC problems may be approximated by Linear Programming (LP) problems without significant increase in the optimization problem complexity. However, $l2$-norm MMMPCs can be represented only by explicit enumeration of all vertexes of the uncertainty polytopic set \cite{gao2012explicit, bemporad2003min}. The enumeration approach suffers from Bellman's curse of dimensionality, since the number of vertexes grows exponentially with the number of dimensions. This issue is further aggravated if the uncertainty set is defined by $l2$-norm, as the enumeration approach cannot be performed on non-polytopic sets.

Xie and Soh \cite{xie1993control} designed the Guaranteed Cost Control (GCC) to address the poor performance of the Linear Quadratic Regulator (LQR) when a Linear Time Invariant (LTI) system model is subject to $l2$-norm bounded parametric uncertainties. It provides a quadratically stable closed-form controller that guarantees an upper bound to $l2$-norm cost functionals. However, the GCC design considers only unconstrained control problem.

GCCs are obtained by either Linear Matrix Inequalities (LMI) \cite{petersen1998optimal} or Algebraic Riccati Equations (ARE) \cite{xie1993control} and the concept was extended to polytopic-bounded $l2$-norm RMPCs problems based on LMIs \cite{quelho2013regions, rosinova2003necessary}. However, a semi-definite programming (SDP) problem must be solved for the guaranteed cost constrained optimization, which increases computational requirements in comparison to LPs, QPs and QCQPs.

In this paper we propose a Guaranteed Cost approach for Robust Model Predictive Control (GCMPC) of linear time invariant systems subject to bounded parametric uncertainties. This approach provides a QCQP-based solution to $l2$-norm Robust Model Predictive Controllers subject to $l2$-norm bounded uncertainties and polytopic constraints and ensures sufficient conditions for stability, feasibility and performance robustness. 

This paper is organized as follows: Section \ref{sec_problem_statement} presents the robust constrained optimal control problem; Section \ref{sec_gsmpc} derives the Guaranteed Cost Model Predictive Controller; Section \ref{sec_example} provides a numerical example; finally, Section \ref{sec_conclusion} addresses the final remarks.

\section{PROBLEM STATEMENT AND PRELIMINARIES}
\label{sec_problem_statement}

This study investigates the formulation of the Guaranteed Cost Model Predictive Control (GCMPC) for linear systems subject to parametric uncertainties. A deterministic LTI system model (Definition \ref{def_clti_system}) and a parametrically uncertain one (Definition \ref{def_ulti_system}) were considered for the formulation.

\begin{definition}
Let the discrete state space LTI model be%
\begin{equation}
x_{k+1} = F x_k + G u_k
\end{equation}%
where $ x_k \in \mathds{R}^n $ is the system state, $ u_k \in \mathds{R}^m $ is the control input, $ F \in \mathds{R}^{n \times n} $ is the state matrix and $ G \in \mathds{R}^{n \times m} $ is the input matrix.
\label{def_clti_system}
\end{definition}

\begin{definition}
\cite{terra2014optimal} 
Let the discrete state space LTI model subject to parametric uncertainties be%
\begin{equation}
\bar{x}_{k+1} = (F + \delta F_k) \bar{x}_k + (G + \delta G_k) u_k
\end{equation}%
where $ \bar{x}_k \in \mathds{R}^n $ is the uncertain system state and $ \delta F_k $ and $ \delta G_k $ are, respectively, the state and input multiplicative uncertainty matrices, such that%
\begin{equation}
\begin{bmatrix}
\delta F_k \; \delta G_k
\end{bmatrix}
= H \Delta_k
\begin{bmatrix}
E_1 \; E_2
\end{bmatrix}
\end{equation}%
with $ \mathcal{W} = \{\delta \mid \delta \in \mathds{R}^{p \times l}, || \delta ||_2 \le 1 \} $, $ H \in \mathds{R}^{n \times p}, H \neq 0 $, $ E_1 \in \mathds{R}^{l \times n} $ and $ E_2 \in \mathds{R}^{l \times m} $.
\label{def_ulti_system}
\end{definition}

Optimality is defined by a quadratic cost functional for LQR and $l2$-norm MPC approaches, given by%
\begin{equation}
J_i(x_i, \textbf{u}, N) = x_N^T P_n x_N + \underset{k = i}{\overset{N-1}{\sum}} x_k^T Q x_k + u_k^T R u_k
\end{equation}
where $ P_N \succeq 0 $, $ Q \succeq 0 $ and $ R \succ 0 $ are symmetric weight matrices, $ \mathbf{u} = \{u_k \mid k \in \mathds{Z}_{\ge 0} \cap [0, N-1]\} $ and the system model satisfies either Definition \ref{def_clti_system}, or \ref{def_ulti_system}.

\begin{definition}
\cite{kalman1960contributions}
Consider the deterministic system model from Definition \ref{def_clti_system}. The Linear Quadratic Regulator is a feedback controller given by the optimal control problem%
\begin{equation}
\begin{matrix*}[l]
J_i^*(x_i, N) = & \underset{\mathbf{u}}{\inf} & J_i(x_0, \mathbf{u}, N)\\
& s.t. & x_{k+1} = F x_k + G u_k
\end{matrix*}
\end{equation}%
which has a closed form solution $ J_i^*(x_0) = x_i^T P_i x_i $ achieved at $ u_k = - K_k x_k $, where%
\begin{equation}
\begin{matrix*}[l]
K_k &= (R + G^T P_{k+1} G)^{-1} G^T P_{k+1} F\\
P_k &= F^T P_{k+1} F + Q - K_k^T (R + G^T P_{k+1} G) K_k \\
\end{matrix*}
\end{equation}
\label{def_lqr}
\end{definition}

\begin{definition}
\cite{camacho2013model}
Consider the deterministic system model from Definition \ref{def_clti_system}. A Model Predictive Control is a feedback controller given by the constrained optimal control problem%
\begin{equation}
\begin{matrix*}[l]
J_i^*(x_i, N) = & \underset{\mathbf{u}}{\inf} & J_i(x_0, \mathbf{u}, N)\\
& s.t. & x_{k+1} = F x_k + G u_k\\
& & [x_k^T, u_k^T]^T \in \mathcal{C}_k
\end{matrix*}
\end{equation}%
where $ \mathcal{C}_k $ defines the feasible set of the system states and control inputs at timestep $ k $. 
\label{def_mpc}
\end{definition}

MPC has no general closed form solution. Therefore, it requires either an online evaluation with a Quadratic Programming solver \cite{lofberg2004yalmip}, or an offline evaluation based on explicit solutions \cite{alessio2009survey}.

\begin{assumption}
The feasible set of variables $ x_k $ and $ u_k $ at a given timestep $ k $, given by $ \left[ x_k^T, u_k^T \right]^T \in \mathcal{C}_k $, is a convex non-empty polytope described by%
\begin{equation}
\mathcal{C}_k = \left\{ z \mid z \in \mathds{R}^{n + m}, C_k z + c_k \le 0 \right\}
\end{equation}%
where $ C_k = [A_k, B_k] $ with $ A_k \in \mathds{R}^{q \times n} $ and $ B_k \in \mathds{R}^{q \times m} $
\label{ass_linear_constraints}
\end{assumption}

Assumption \ref{ass_linear_constraints} enables the use of widely available QP solvers (e.g. Gurobi \cite{gurobi}, CVXGEN \cite{mattingley2012cvxgen} and FORCES Pro \cite{FORCESPro}) and the solution of several practical control problems. Therefore, it is used in most MPC applications. However, neither LQR, nor MPC provides performance, stability or feasibility robustness guarantees when they are used for the control of systems subject to multiplicative or additive uncertainties.

\begin{definition}
\cite{petersen1998optimal}
Consider the system model from Definition \ref{def_ulti_system}. A control law $ u_k = - K_k \bar{x}_k $ is said to be a Guaranteed Cost Controller if there exists a symmetric cost matrix sequence $ \mathcal{S} = \{S_k \mid S_k \succeq 0, k \in \mathds{Z}, k \in [0, N]\} $, such that%
\begin{multline}
\left[ F - G K_k + \delta F_k - \delta G_k K_k \right]^T S_{k+1} \left[ \bullet \right] - \\ - S_k + Q + K_k^T R K_k < 0
\label{eq_guarateeed_cost_lyapunov}
\end{multline}%
where $ z^T M ( \bullet ) = z^T M z $.
\label{def_guaranteed_cost}
\end{definition}

By definition, GCC ensures the quadratic stability of the system for any disturbance within the admissible set $ \mathcal{W} $. A closed form solution for $ \mathcal{S} $ based on AREs is proposed by Xie and Soh \cite{xie1993control}.

\begin{theorem}
Consider the parametrically uncertain system model from Definition \ref{def_ulti_system} and assume a control law $ u_k = - K_k \bar{x}_k $ is a Guaranteed Cost Control with an associated cost matrix $ S_k $. For some $ \epsilon_k > 0 $, such that $ S_k \succ 0 $, the closed loop uncertain system%
\begin{equation}
\bar{x}_{k+1} = [F - G K_k + \delta F_k - \delta G_k K_k] \bar{x}_k
\end{equation}%
has an upper bound cost%
\begin{equation}
J_k(\bar{x}_k, \mathbf{u}, N) < \bar{x}_k^T S_k \bar{x}_k
\end{equation}%
given by%
\begin{multline}
S_k = F^T X_{k+1} F + Q + \epsilon_k^{-1} E_1^T E_1 - \\ - (F^T X_{k+1} G + \epsilon_k^{-1} E_1^T E_2) (R + \epsilon_k^{-1} E_2^T E_2 + G^T X_{k+1} G)^{-1} (\bullet)
\end{multline}%
where $ I - \epsilon_k H^T S_k H \succ 0 $ and
\begin{align}
X_k =& \left( S_k^{-1} - \epsilon_k^{-1} H H^T\right)^{-1} \succ 0\\
K_k =& (R_{\epsilon k} + G^T X_{k+1} G)^{-1} (G^T X_{k+1} F + \epsilon_k^{-1} E_2^T E_1).
\end{align}
\label{the_guaranteed_cost}

\begin{proof}
See \cite{xie1993control}.
\end{proof}

\end{theorem}

However, the results from Theorem \ref{the_guaranteed_cost} are valid only for $ [\bar{x}_k^T, u_k^T]^T \in \mathds{R}^{n+m} $ and have not been extended to the constrained quadratic optimal control case.

\section{GUARANTEED COST MODEL PREDICTIVE CONTROL}
\label{sec_gsmpc}

This section addresses the formulation of the constrained guaranteed cost control.

\begin{lemma}
Let $ \epsilon \in \left(0, ||H^T S H|| \right) $ and%
\begin{equation}
X = \left( S^{-1} - \epsilon H H^T \right)^{-1} \succ 0
\label{eq_x_s_relation}
\end{equation}%
for a given symmetric $ S \succ 0 $ and a non-singular $ H $. Then, $X \succ S$.

\begin{proof}
We employ Woodbury's matrix inversion Lemma \cite{henderson1981deriving} to \eqref{eq_x_s_relation}, which results in%
\begin{equation}
X = S + S H \left( \epsilon^{-1} I - H^T S H \right)^{-1} H^T S
\end{equation}%
and conclude, from $ \epsilon^{-1} I - H^T S H \succ 0 $,
\begin{equation}
S H \left( \epsilon^{-1} I - H^T S H \right)^{-1} H^T S \succ 0
\end{equation}%
Therefore, $ X \succ S $.
\end{proof}
\label{lem_x_gt_s}
\end{lemma}

\begin{theorem}
Consider the deterministic system model from Definition \ref{def_clti_system} and%
\begin{equation}
\begin{matrix*}[l]
\bar{J}_i^*(x_i, N) = & \underset{\mathbf{u}}{\inf} & \bar{J}_i(x_i, \mathbf{u}, N)\\
& s.t. & x_{k+1} = F x_k + G u_k
\label{eq_opt_gcc}
\end{matrix*}
\end{equation}%
where%
\begin{align}
\bar{J}_i(x_i, \mathbf{u}, N) &= x_N^T S_N x_N + \underset{k = i}{\overset{N}{\sum}} c(x_k, u_k) \\
c(x, u) &= x^T Q_{\epsilon k} x + u^T R_{\epsilon_k} u + 2 x^T N_{\epsilon_k} u
\end{align}%
with $ Q_{\epsilon k} = Q + \epsilon^{-1} E_1^T E_1 $, $ R_{\epsilon k} = R + \epsilon^{-1} E_2^T E_2 $, $ N_{\epsilon k} = \epsilon^{-1} E_1^T E_2 $ and $ \epsilon $ is chosen such that the guaranteed cost control condition is satisfied. The guaranteed cost control is equivalent to the sub-optimal solution of \eqref{eq_opt_gcc}, given by%
\begin{equation}
\begin{matrix*}[l]
\bar{J}_N(x_N) &< x_N^T X_N x_N\\
\bar{J}_k(x_k) &< \underset{u_k}{\inf} \; c(x_k, u_k) + x_{k+1}^T X_{k+1} x_{k+1} = x_k^T S_k x_k
\end{matrix*}
\label{eq_value_function_gcc}
\end{equation}

\begin{proof}
Let the value function $ \bar{J}^*_k(x) $, related to the optimization problem from \eqref{eq_opt_gcc}, be recursively defined as %
\begin{align}
\label{eq_J_N}
\bar{J}^*_N(x) &= x^T S_N x \\
\bar{J}^*_k(x) &= \underset{u}{\inf} \; c(x, u) + \bar{J}_{k+1}(F x + G u)
\label{eq_J_k}
\end{align}%
based on Bellman's optimality principle.

From Lemma \ref{lem_x_gt_s} and \eqref{eq_J_N},%
\begin{equation}
\bar{J}^*_N(x) = x^T S_N x < x^T X_N x
\end{equation}%
Therefore,%
\begin{align}
\bar{J}^*_{N-1}(x) &= \underset{u}{\inf} \; c(x, u) + (F x + G u)^T S_{N} ( \bullet ) \\
&< \underset{u}{\inf} \; c(x, u) + (F x + G u)^T X_{N} ( \bullet )
\label{eq_J_k_ineq}
\end{align}

The solution of \eqref{eq_J_k_ineq} yields%
\begin{equation}
u = - (R_{\epsilon_k} + G^T X_N G)^{-1} (G^T X_N F + N_{\epsilon_k}) x
\end{equation}%
which substituted back in \eqref{eq_J_k_ineq} results in%
\begin{multline}
\bar{J}^*_{N-1}(x) < x^T \left[ F^T X_{N} F + Q_{\epsilon k} - \right. \\
\left. - (F^T X_{N} G + N_{\epsilon k}) (R_{\epsilon k} + G^T X_{N} G)^{-1} (\bullet) \right] x
\end{multline}%
Therefore, $ \bar{J}^*_{N-1}(x) < x^T S_{N-1} x $, where $ S_{N-1} $ is the guaranteed cost matrix at timestep $ k = N-1 $, as stated in Theorem \ref{the_guaranteed_cost}.

Assume $ \exists k \in \mathds{Z}_{\ge 0} \cap [0, N-2]: \; \bar{J}^*_{k+1}(x) < x^T S_{k+1} x $, where $ S_{k+1} $ is a guaranteed cost matrix.%
\begin{align}
\bar{J}^*_{k}(x) &< \underset{u}{\inf} \; c(x, u) + (F x + G u)^T S_{k+1} ( \bullet ) \\
&< \underset{u}{\inf} \; c(x, u) + (F x + G u)^T X_{k+1} ( \bullet )
\end{align}%
which, analogously to $ \bar{J}_{N-1} $, results in%
\begin{multline}
\bar{J}_{k}(x) < x^T \left[ F^T X_{k+1} F + Q_{\epsilon k} - \right. \\
\left. - (F^T X_{k+1} G + N_{\epsilon k}) (R_{\epsilon k} + G^T X_{k+1} G)^{-1} (\bullet) \right] x
\end{multline}%
Therefore, we conclude by induction%
\begin{equation}
\forall k \in \mathds{Z}_{\ge 0} \cap [0, N-1]. \; \bar{J}_k(x) < x^T S_k x
\end{equation} %
hence the guaranteed cost control is a sub-optimal solution to the optimization problem in \eqref{eq_opt_gcc}.
\end{proof}
\label{the_subopt_gcc}
\end{theorem}

Theorem \ref{the_subopt_gcc} demonstrates the equivalence of the GCC of a parametrically uncertain system to a sub-optimal solution of a Dynamic Programming problem for the deterministic system model. It also shows that the GCC bounds $ J_i^*(\bar{x}_i, \mathbf{u}, N) < x_i^T S_i x_i $ and $ \bar{J}_i^*(x_i, \mathbf{u}, N) < x_i^T S_i x_i $. However, the result from Theorem \ref{the_subopt_gcc} is valid only for the unconstrained case where $ x_k \in \mathds{R}^n $ and $ u_k \in \mathds{R}^m $.

\begin{lemma}
Consider the parametrically uncertain system model from Definition \ref{def_ulti_system} and let $ K_k $ be a guaranteed cost control, $ u_k = - K_k \bar{x}_k + v_k $ and $ \bar{R}_k = R_{\epsilon k} + G^T X_{k+1} G $.%
\begin{equation}
\begin{matrix*}[l]
J_i^*(x_i, \mathbf{u}, N) < & \underset{\mathbf{v}}{\inf} & \underset{k = i}{\overset{N-1}{\sum}} v_k \bar{R}_k v_k + x_i^T S_i x_i\\
& s.t. & x_{k+1} = (F - G K_k) x_k + G v_k\\
\end{matrix*}
\label{eq_gsmpc_cost}
\end{equation}%
is equivalent to the guaranteed cost control sub-optimal solution from Theorem \ref{the_subopt_gcc}.

\begin{proof}
Consider the value function from \eqref{eq_value_function_gcc} and the deterministic system model from Definition \ref{def_clti_system}. %
\begin{multline}
\bar{J}^*_k(x) < \underset{u}{\inf} \; x^T (Q_{\epsilon k} + F^T X_{k+1} F) x + \\ + u^T (R_{\epsilon k} + G^T X_{k+1} G) u + 2 x^T F^T X_{k+1} G u
\end{multline}%
where the substitution of $ u = K_k x_k + v_k $ yields%
\begin{equation}
\bar{J}^*_k(x) < x^T S_k x + \underset{v}{\inf} \; v^T (R_{\epsilon k} + G^T X_{k+1} G) v
\label{eq_min_subopt_gcc_value}
\end{equation}

The recursive substitution of \eqref{eq_min_subopt_gcc_value} for $ k \in [i, N-1] $ results in%
\begin{equation}
\begin{matrix*}[l]
\bar{J}^*_i(x_i) < & \underset{\mathbf{v}}{\inf} & \underset{k = i}{\overset{N-1}{\sum}} v_k \bar{R} v_i + x_i^T S_i x_i\\
& s.t. & x_{k+1} = (F - G K_k) x_k + G v_k\\
\end{matrix*}
\end{equation}%
Therefore, from the results of Theorem \ref{the_subopt_gcc} and \cite{xie1993control},%
\begin{equation}
\begin{matrix*}[l]
J_i^*(x_i, \mathbf{u}, N) < & \underset{\mathbf{v}}{\inf} & \underset{k = i}{\overset{N-1}{\sum}} v_k \bar{R} v_i + x_i^T S_i x_i\\
& s.t. & x_{k+1} = (F - G K_k) x_k + G v_k\\
\end{matrix*}
\label{eq_unconstrained_gcmpc}
\end{equation}
\end{proof}
\label{the_gsmpc_cost}
\end{lemma}

If the optimization domain of \eqref{eq_gsmpc_cost} is $ \mathcal{C}_k \equiv \mathds{R}^{n+m} $, the solution is trivially $ v_k \equiv 0 $ and \eqref{eq_gsmpc_cost} reduces to Theorem \ref{the_subopt_gcc}. However, if $ \mathcal{C}_k \subset \mathds{R}^{n+m} $, \eqref{eq_gsmpc_cost} provides an upper bound to the MPC cost functional.

Since \eqref{eq_unconstrained_gcmpc} is based on the deterministic model, enforcement of the feasible set constraints ($ [x_k^T, u_k^T]^T \in \mathcal{C}_k $) would ensure that only the undisturbed states and control inputs are feasible and there would be no feasibility guarantees for the parametrically uncertain ones. Feasibility robustness requires%
\begin{equation}
\forall \Delta_k \in \mathcal{W}. \; [\bar{x}_k^T, u_k^T]^T \in \mathcal{C}_k
\label{eq_uncertain_feasible_set}
\end{equation}

\begin{lemma}
Consider the parametrically uncertain system model from Definition \ref{def_ulti_system} and let $ K_k $ be a Guaranteed Cost Controller, $ u_k = - K_k \bar{x}_k + v_k $ and $ \mathcal{\bar{W}}_k = \{w \mid w \in \mathds{R}^{p \times l}, || w ||_2 \le ||\widetilde{E}_1 \bar{x}_k + E_2 v_k||_2 \} $ where $ \widetilde{E}_1 = E_1 - E_2 K_k $.%
\begin{equation}
\bar{x}_{k+1} = (F - G K_k) \bar{x}_k + G v_k + H w_k
\label{eq_ulti_wk}
\end{equation}%
is an approximation of the closed loop parametrically uncertain model.

\begin{proof}
Let $ w_k = \Delta_k [\widetilde{E}_{1,k} \bar{x}_k + E_2 v_k] $. From \eqref{eq_ulti_wk}, we conclude%
\begin{equation}
\begin{matrix*}[l]
\bar{x}_{k+1} &= (F - G K_k) \bar{x}_k + G v_k + H \Delta_k (\widetilde{E}_{1,k} \bar{x}_k + E_2 v_k) \\
&= F \bar{x}_k + G u_k + H \Delta_k (E_1 \bar{x}_k + E_2 u_k)  \\
&= (F + \delta F_k) \bar{x}_k + (G + \delta G_k) u_k
\end{matrix*}
\end{equation}%
and the admissible set of $ w_k $, $ \mathcal{\bar{W}} $, is given by%
\begin{equation}
\begin{matrix*}[l]
|| w_k ||_2 &= || \Delta_k [\widetilde{E}_{1,k} \bar{x}_k + E_2 v_k] ||_2 \\
& \le ||\Delta_k||_2 ||\widetilde{E}_{1,k} \bar{x}_k + E_2 v_k||_2 \\
& \le ||\widetilde{E}_{1,k} \bar{x}_k + E_2 v_k||_2
\end{matrix*}
\end{equation}%
Therefore, $ \mathcal{\bar{W}}_k = \{w \mid w \in \mathds{R}^{p \times l}, || w ||_2 \le ||\widetilde{E}_1 \bar{xs}_k + E_2 v_k||_2 \} $. Since $ \mathcal{\bar{W}} \supset \mathcal{W} $, such a representation is a conservative approximation of Definition \ref{def_ulti_system}.
\end{proof}
\label{lem_ulti_approx}
\end{lemma}

Lemma \ref{lem_ulti_approx} provides a conservative approximation to the parametrically uncertain system. Such representation enables multiplicative uncertainties to be modeled as in previous RMPC studies (i.e. \cite{scokaert1998min}).

\begin{corollary}
Consider the deterministic and parametrically uncertain system models from Definition \ref{def_clti_system} and Lemma \ref{lem_ulti_approx}, respectively, and assume the system is regulated by a feedback controller $ \widetilde{K} $.%
\begin{equation}
\bar{x_k} = x_k + \underset{i = 0}{\overset{k - 1}{\sum}} \widetilde{F}^{k - i - 1} H w_i
\end{equation}%
where $ \widetilde{F} = F - G \widetilde{K} $ is the closed loop system dynamics.
\label{cor_system_relation}
\end{corollary}

Corollary \ref{cor_system_relation} provides a representation of the uncertain system state according to the deterministic system state and all previous disturbances and enables the externalization of disturbances from the system dynamics to the feasible set constraints. Therefore, based on \eqref{eq_uncertain_feasible_set} and Assumption \ref{ass_linear_constraints}, the parametrically uncertain model feasible set $ \mathcal{C}_k $ with $ u_k = - K_k x_k + v_k $ is%
\begin{multline}
\forall w_k \in \mathcal{\bar{W}}. \; \widetilde{A}_k x_k + B_k v_k + c_k + \\ + \widetilde{A}_k \underset{i = 0}{\overset{k - 1}{\sum}} \widetilde{F}^{k - i - 1} H w_i \le 0
\label{eq_uncertain_constraint_replaced}
\end{multline}%
where $ \widetilde{A}_k = A_k - B_k K_k $ and $ \widetilde{K} $ is arbitrarily defined. However, uncertainty admissible set $ \mathcal{\bar{W}} $ is still dependent on uncertain state $ \bar{x}_k $.

\begin{lemma}
Let $ \phi_k(x, v) = || \widetilde{E}_{1,k} x + E_2 v ||_2 $ and $ \rho_i = || E_1 \widetilde{F}^i H ||_2 $. Then,%
\begin{multline}
|| w_k || \le \phi_k(\bar{x}_k, v_k) \le \\ \phi_k(x_k, v_k) + \underset{i = 0}{\overset{k - 1}{\sum}} c(k,i) \phi_i(x_i, v_i)
\end{multline}%
where%
\begin{equation}
\forall i < k. \; c(k,i) = \rho_{k-i-1} + \underset{j = 0}{\overset{k - i - 2}{\sum}} \rho_j c(k-j-1, i) 
\end{equation}

\begin{proof}
Omitted.
\end{proof}
\label{lem_phi_function}
\end{lemma}

Duality filter, enumeration filter, explicit maximization filter, Pólya filter and elimination filter are five main approaches that devise robust counterparts of $ \mathcal{C}_k $ \cite{lofberg2012automatic}, such that \eqref{eq_uncertain_feasible_set} is satisfied for an optimization based on the deterministic model. Pólya and elimination filters are relaxations of the constraints and are not sufficient conditions to satisfy \eqref{eq_uncertain_constraint_replaced}. Enumeration and the duality filters cause the number of variables to grow significantly. Therefore, this study has focused on the explicit maximization filter.

\begin{definition}
Given an uncertain variable $ w \in \{w \mid w \in \mathds{R}^{n}, || w ||_p \le 1\} $ and an inequality constraint%
\begin{equation}
(A^T w + b)^T x + (c^T w + d) \le 0
\end{equation}%
The robust constraint counterpart by explicit maximization is given by%
\begin{equation}
b^T x + d + || A x + c ||_{p^*} \le 0
\end{equation}%
where $ || z ||_{p^*} $ is the dual norm \cite{boyd2004convex} of $ || A x + c ||_p $ and satisfies $ 1 / p + 1 / p^* = 1 $.
\label{def_explicit_maximization}
\end{definition}

\begin{theorem}
Let $ \bar{\mathcal{C}}_k $ be the robust counterpart of $ \mathcal{C}_k $ by explicit maximization and, based on Lemma \ref{lem_phi_function}, $ \bar{\phi}_k(\mathbf{x}, \mathbf{v}) = \phi_k(x_k, v_k) + \underset{i = 0}{\overset{k - 1}{\sum}} c(k,i) \phi_i(x_i, v_i) $. Then,%
\begin{multline}
\bar{\mathcal{C}}_k = \{ [x_k^T, v_k^T]^T \mid x_k \in \mathds{R}^n, v_k \in \mathds{R}^m, \forall i \in \mathds{I} \cap [1, q]. \\ \widetilde{A}_k^{(i)} x_k + B_k^{(i)} v_k + c_k^{(i)} + \Phi_{k,i}(\mathbf{x}, \mathbf{v}) \le 0 \}
\end{multline}%
where $ \Phi_{k,i}(\mathbf{x}, \mathbf{v}) = \underset{j = 0}{\overset{k - 1}{\sum}} || \widetilde{A}_k^{(i)} \widetilde{F}^{k - j - 1} H ||_2 \bar{\phi}_j(\mathbf{x}, \mathbf{v}) $ and $ \widetilde{A}_k^{(i)} $, $ B_k^{(i)} $ and $ c_k^{(i)} $ are the $i$-th row of matrices $ \widetilde{A}_k $, $B_k$ and vector $c_k$, respectively.

\begin{proof}
Consider Definition \ref{def_explicit_maximization} and Lemma \ref{lem_phi_function}. The explicitly maximized counterpart of $ \widetilde{A}_k^{(i)} \widetilde{F}^{k - i - 1} H w_i $ is%
\begin{multline}
\underset{||w_i||_2 \le \bar{\phi}_k(\mathbf{x}, \mathbf{v})}{\max} \widetilde{A}_k^{(i)} \widetilde{F}^{k - i - 1} H w_i = \\ = || \widetilde{A}_k^{(i)} \widetilde{F}^{k - i - 1} H ||_2 \bar{\phi}_k(\mathbf{x}, \mathbf{v})
\label{eq_explicit_maximization_proof_1}
\end{multline}

Therefore, the explicit maximization of \eqref{eq_uncertain_constraint_replaced} results%
\begin{multline}
\widetilde{A}_k^{(i)} x_k + B_k^{(i)} v_k + c_k^{(i)} + \\ + \underset{j = 0}{\overset{k - 1}{\sum}}  || \widetilde{A}_k^{(i)} \widetilde{F}^{k - j - 1} H ||_2 \bar{\phi}_j(\mathbf{x}, \mathbf{v}) \le 0
\end{multline}%
which is a sufficient condition for the robust feasibility of \eqref{eq_uncertain_constraint_replaced}.

\end{proof}
\label{the_robust_feasible_set}
\end{theorem}

The results from Theorems \ref{the_subopt_gcc} and \ref{the_robust_feasible_set} and Lemma \ref{the_gsmpc_cost} provide the basis of the proposed method.

\begin{corollary}
Consider the results from Theorem \ref{the_robust_feasible_set} and Lemma \ref{the_gsmpc_cost}. The guaranteed cost model predictive controller is given by%
\begin{equation}
\begin{matrix*}[l]
J_i^*(\bar{x}_i, N) < & \underset{\mathbf{v}}{\inf} & \underset{k = i}{\overset{N-1}{\sum}} v_k \bar{R}_k v_k + x_i^T S_i x_i\\
& s.t. & x_{k+1} = (F - G K_k) x_k + G v_k\\
& & [x_k^T, v_k^T]^T \in \bar{\mathcal{C}}_k
\end{matrix*}
\end{equation}
\label{cor_gsmpc}
\end{corollary}

Corollary \ref{cor_gsmpc} defines the Guaranteed Cost Model Predictive Control, which ensures an upper bound to the linearly constrained $l2$-norm optimal control problem subject to parametric uncertainties. 

\section{NUMERICAL EXAMPLE}
\label{sec_example}

This section provides a numerical example of GCMPC and its comparison with a Enumeration-based RMPC (ERMPC) approach\footnotemark[2], adapted from \cite{scokaert1998min}. The YALMIP Toolbox \cite{lofberg2004yalmip} and the Mosek solver \cite{mosek} were used in the modeling and simulation of the problem\footnotemark[1].

\footnotetext[1]{The source code for the numerical example with both GCMPC and ERMPC  are available at: \url{https://github.com/cmasseraf/gcmpc}}
\footnotetext[2]{The approximation $ ||\Delta||_\infty \le ||\Delta||_2 \le 1 $ was used to transform the disturbance set into a polytope and enable enumeration to be performed.} Consider the system from Definition \ref{def_ulti_system} with matrices%
\begin{equation}
\begin{matrix}
\begin{matrix}
F = \begin{bmatrix}
1.1 & 0 & 0 \\
0 & 0 & 1.2 \\
-1 & 1 & 0
\end{bmatrix} & 
G = \begin{bmatrix}
0 & 1 \\
1 & 1 \\
-1 & 0
\end{bmatrix}
\end{matrix}\\
\begin{matrix}
H = \begin{bmatrix}
0.7 \\ 0.5 \\ -0.7
\end{bmatrix} &
E1 = \begin{bmatrix}
0.4 \\ 0.5 \\ -0.6
\end{bmatrix}^T & 
E2 = \begin{bmatrix}
0.4 \\ -0.4
\end{bmatrix}^T
\end{matrix}
\end{matrix}
\end{equation}%
subject to constraints given by%
\begin{equation}
\begin{matrix}
A_k \equiv \begin{bmatrix}
1 & 0 & 0\\
-1 & 0 & 0\\
0 & 1 & 0\\
0 & -1 & 0\\
0 & 0 & 1\\
0 & 0 & -1\\
\end{bmatrix} &
B_k \equiv \mathds{O}_6 & 
c_k \equiv - \mathds{1}_6
\end{matrix}
\end{equation}%
where $ \mathds{O}_i $ and $ \mathds{1}_i $ are zero and one valued column vectors of size $ i $, respectively, and $ Q = \mathds{I}_3 $, $ R = \mathds{I}_2 $ and $ N = 10 $.

\begin{figure}[!t]
  \centering
  \includegraphics[width=0.9\columnwidth]{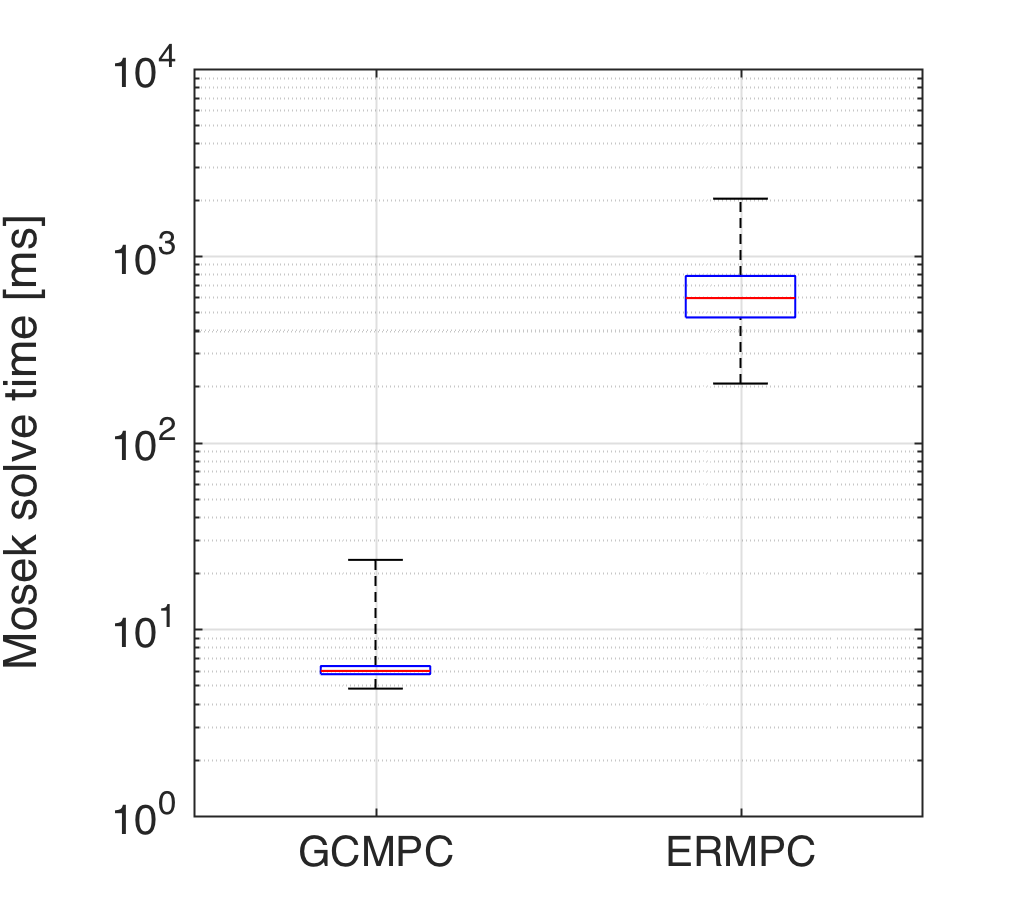}
  \caption{Computational time box plot of both GCMPC and ERMPC for uniformly distributed states with $||x||_\infty \le 0.5 $ }
  \label{fig_comp_time}
\end{figure}

\begin{figure*}[!t]
\centering
\begin{subfigure}{0.3\textwidth}
  \centering
  \includegraphics[width=\linewidth]{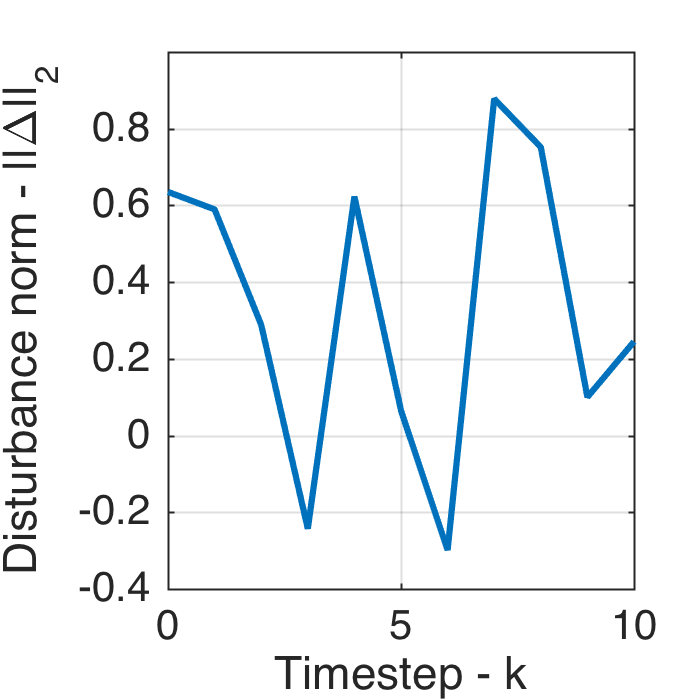}
  \caption{}
  \label{fig_result_a}
\end{subfigure}
\quad
\begin{subfigure}{0.3\textwidth}
  \centering
  \includegraphics[width=\linewidth]{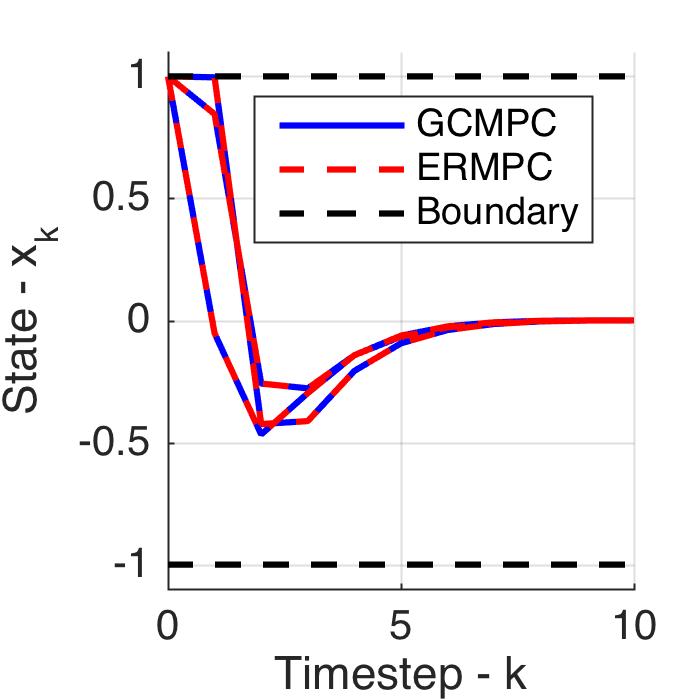}
  \caption{}
  \label{fig_result_b}
\end{subfigure}
\quad
\begin{subfigure}{0.3\textwidth}
  \centering
  \includegraphics[width=\linewidth]{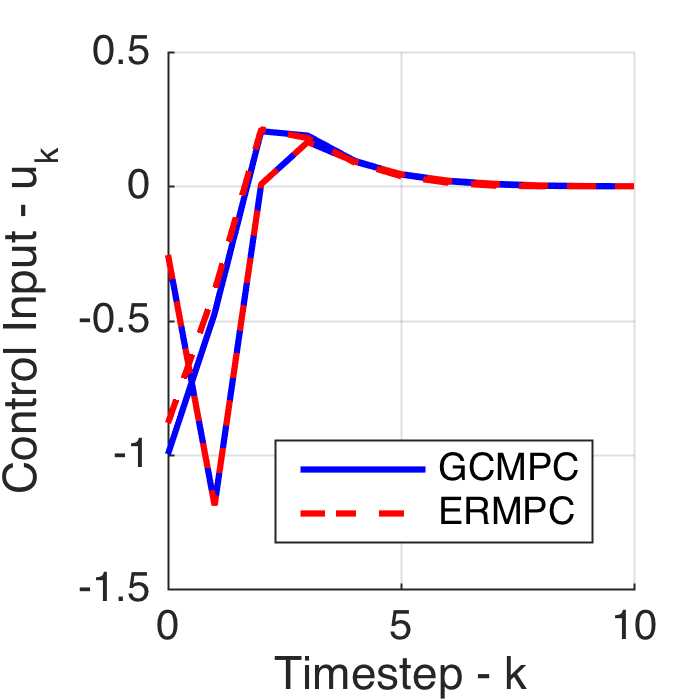}
  \caption{}
  \label{fig_result_c}
\end{subfigure}
\caption{Simulation results for GCMPC and ERMPC controllers. (a) Simulated disturbance norm $ || \Delta ||_2 $; (b) Disturbed closed-loop state $ x_k $; (c) Executed control inputs $ u_k $.}
\label{fig_result}
\end{figure*}

An infinite horizon positive definite solution $ S $ exists to the GCC problem for $ \epsilon_k \equiv \epsilon \in (0, 0.0220) $ and $ \epsilon = 0.0180 $ provides the optimal stabilizing solution, given by%
\begin{equation}
S_k \equiv S = \begin{bmatrix}
31.4751 & -0.9359 & -20.6124 \\
-0.9359 & 5.7340 & -1.3900 \\
-20.6124 & -1.3900 & 16.5017
\end{bmatrix}
\end{equation}%
and%
\begin{equation}
K_k \equiv K = \begin{bmatrix}
1.1801 & 0.2151 & -0.5076 \\
0.7401 & -0.8385 & 0.5162
\end{bmatrix}
\end{equation}%
whereas the GCMPC cost matrix $ \bar{R}_k $ is%
\begin{equation}
\bar{R}_k \equiv \bar{R} = \begin{bmatrix}
123.22 & 133.78 \\
133.78 & 197.26
\end{bmatrix}
\end{equation}%
and $ \widetilde{K} $, a second order nil-potent controller, is%
\begin{equation}
\widetilde{K} = \begin{bmatrix}
11.5 & -6 & -6 \\
1.1 & 0 & 0
\end{bmatrix}
\end{equation}

The resulting optimization problems were solved in a minimum of $4.8ms$ and a maximum of $23.5ms$ for the GCMPC and a minimum of $ 206.77 ms$ and a maximum of $2021.5 ms$ for the ERMPC on a $2.5Ghz$ i7-4980HQ with $16Gb$ of RAM. The computational time distributions is shown in Figure \ref{fig_comp_time}, where it is possible to see that the GCMPC is on average two orders of magnitude faster than the ERMPC.

Figure \ref{fig_result} shows the simulation results of both controllers for an initial state $ x_0 = [1, 1, 1]^T $. The controlled system was subject to a uniformly distributed disturbance $ \Delta_k \in [-1, 1] $, shown in Figure \ref{fig_result_a}). The GCMPC successfully maintained the system bounded and stabilized with the an equivalent  performance of the optimal ERMPC, as shown in Figure \ref{fig_result_b}.

\section{CONCLUSIONS}
\label{sec_conclusion}

In this paper we have proposed a constrained guaranteed cost robust model predictive controller for uncertain discrete time systems.
Such controller provides an upper bound to the quadratic cost functional and ensures feasibility of states and control inputs within a prediction horizon for a system subject to a quadratically bound multiplicative uncertainty. This problem class was only tractable though vertex enumeration of a polytopic approximation of the uncertainty set.

The proposed approach is computationally more efficient than enumeration-based techniques and executed on average two orders of magnitude faster for this paper numerical example. It also requires only a QCQP optimization, which results in lower-complexity solutions than those of LMI-based guaranteed cost approaches based on SDP optimizations.


\bibliographystyle{IEEEtranS}
\bibliography{refs_paper}

\begin{thebibliography}{10}
\providecommand{\url}[1]{#1}
\csname url@samestyle\endcsname
\providecommand{\newblock}{\relax}
\providecommand{\bibinfo}[2]{#2}
\providecommand{\BIBentrySTDinterwordspacing}{\spaceskip=0pt\relax}
\providecommand{\BIBentryALTinterwordstretchfactor}{4}
\providecommand{\BIBentryALTinterwordspacing}{\spaceskip=\fontdimen2\font plus
\BIBentryALTinterwordstretchfactor\fontdimen3\font minus
  \fontdimen4\font\relax}
\providecommand{\BIBforeignlanguage}[2]{{%
\expandafter\ifx\csname l@#1\endcsname\relax
\typeout{** WARNING: IEEEtranS.bst: No hyphenation pattern has been}%
\typeout{** loaded for the language `#1'. Using the pattern for}%
\typeout{** the default language instead.}%
\else
\language=\csname l@#1\endcsname
\fi
#2}}
\providecommand{\BIBdecl}{\relax}
\BIBdecl

\bibitem{alessio2009survey}
A.~Alessio and A.~Bemporad, ``A survey on explicit model predictive control,''
  in \emph{Nonlinear model predictive control}.\hskip 1em plus 0.5em minus
  0.4em\relax Springer, 2009, pp. 345--369.

\bibitem{mosek}
\BIBentryALTinterwordspacing
M.~ApS, \emph{The MOSEK optimization toolbox for MATLAB manual. Version 7.1
  (Revision 28).}, 2015. [Online]. Available:
  \url{http://docs.mosek.com/7.1/toolbox/index.html}
\BIBentrySTDinterwordspacing

\bibitem{badgwell2015model}
T.~A. Badgwell and S.~J. Qin, ``Model-predictive control in practice,''
  \emph{Encyclopedia of Systems and Control}, pp. 756--760, 2015.

\bibitem{bemporad2003min}
A.~Bemporad, F.~Borrelli, and M.~Morari, ``Min-max control of constrained
  uncertain discrete-time linear systems,'' \emph{Automatic Control, IEEE
  Transactions on}, vol.~48, no.~9, pp. 1600--1606, 2003.

\bibitem{bemporad1999robust}
A.~Bemporad and M.~Morari, ``Robust model predictive control: A survey,'' in
  \emph{Robustness in identification and control}.\hskip 1em plus 0.5em minus
  0.4em\relax Springer, 1999, pp. 207--226.

\bibitem{boyd2004convex}
S.~Boyd and L.~Vandenberghe, \emph{Convex Optimization}.\hskip 1em plus 0.5em
  minus 0.4em\relax New York, NY, USA: Cambridge University Press, 2004.

\bibitem{camacho2013model}
E.~F. Camacho and C.~B. Alba, \emph{Model predictive control}.\hskip 1em plus
  0.5em minus 0.4em\relax Springer Science \& Business Media, 2013.

\bibitem{gao2012explicit}
Y.~Gao and K.~T. Chong, ``The explicit constrained min-max model predictive
  control of a discrete-time linear system with uncertain disturbances,''
  \emph{Automatic Control, IEEE Transactions on}, vol.~57, no.~9, pp.
  2373--2378, 2012.

\bibitem{genceli1993robust}
H.~Genceli and M.~Nikolaou, ``Robust stability analysis of constrained l1-norm
  model predictive control,'' \emph{AIChE Journal}, vol.~39, no.~12, pp.
  1954--1965, 1993.

\bibitem{FORCESPro}
E.~GmbH, ``{FORCES Professional},''
  {(\nobreak{\url{http://embotech.com/FORCES-Pro}})}, Apr. 2016.

\bibitem{gurobi}
\BIBentryALTinterwordspacing
I.~Gurobi~Optimization, ``Gurobi optimizer reference manual,'' 2015. [Online].
  Available: \url{http://www.gurobi.com}
\BIBentrySTDinterwordspacing

\bibitem{henderson1981deriving}
H.~V. Henderson and S.~R. Searle, ``On deriving the inverse of a sum of
  matrices,'' \emph{Siam Review}, vol.~23, no.~1, pp. 53--60, 1981.

\bibitem{kalman1960contributions}
R.~E. Kalman \emph{et~al.}, ``Contributions to the theory of optimal control,''
  \emph{Bol. Soc. Mat. Mexicana}, vol.~5, no.~2, pp. 102--119, 1960.

\bibitem{kerrigan2004feedback}
E.~C. Kerrigan and J.~M. Maciejowski, ``Feedback min-max model predictive
  control using a single linear program: robust stability and the explicit
  solution,'' \emph{International Journal of Robust and Nonlinear Control},
  vol.~14, no.~4, pp. 395--413, 2004.

\bibitem{lofberg2003approximations}
J.~L{\"o}fberg, ``Approximations of closed-loop minimax {MPC},'' in
  \emph{Decision and Control, 2003. Proceedings. 42nd IEEE Conference on},
  vol.~2.\hskip 1em plus 0.5em minus 0.4em\relax IEEE, 2003, pp. 1438--1442.

\bibitem{lofberg2004yalmip}
------, ``Yalmip: A toolbox for modeling and optimization in matlab,'' in
  \emph{Computer Aided Control Systems Design, 2004 IEEE International
  Symposium on}.\hskip 1em plus 0.5em minus 0.4em\relax IEEE, 2004, pp.
  284--289.

\bibitem{lofberg2012automatic}
------, ``Automatic robust convex programming,'' \emph{Optimization methods and
  software}, vol.~27, no.~1, pp. 115--129, 2012.

\bibitem{mattingley2012cvxgen}
J.~Mattingley and S.~Boyd, ``Cvxgen: A code generator for embedded convex
  optimization,'' \emph{Optimization and Engineering}, vol.~13, no.~1, pp.
  1--27, 2012.

\bibitem{petersen1998optimal}
I.~R. Petersen, D.~C. McFarlane, and M.~A. Rotea, ``Optimal guaranteed cost
  control of discrete-time uncertain linear systems,'' \emph{International
  Journal of Robust and Nonlinear Control}, vol.~8, no.~8, pp. 649--657, 1998.

\bibitem{qin2003survey}
S.~J. Qin and T.~A. Badgwell, ``A survey of industrial model predictive control
  technology,'' \emph{Control engineering practice}, vol.~11, no.~7, pp.
  733--764, 2003.

\bibitem{quelho2013regions}
F.~Quelho~Rossi, R.~Waschburger, and R.~Kawakami Harrop~Galvao, ``Regions of
  guaranteed cost for lmi-based robust model predictive controllers for systems
  with uncertain input delay,'' in \emph{Control and Fault-Tolerant Systems
  (SysTol), 2013 Conference on}.\hskip 1em plus 0.5em minus 0.4em\relax IEEE,
  2013, pp. 584--589.

\bibitem{rawlings1994nonlinear}
J.~Rawlings, E.~Meadows, and K.~Muske, ``Nonlinear model predictive control: A
  tutorial and survey,'' \emph{Advanced Control of Chemical Processes}, pp.
  203--214, 1994.

\bibitem{rosinova2003necessary}
D.~Rosinov{\'a}, V.~Vesel{\`y}, and V.~Ku{\v{c}}era, ``A necessary and
  sufficient condition for static output feedback stabilizability of linear
  discrete-time systems,'' \emph{Kybernetika}, vol.~39, no.~4, pp. 447--459,
  2003.

\bibitem{scokaert1998min}
P.~Scokaert and D.~Mayne, ``Min-max feedback model predictive control for
  constrained linear systems,'' \emph{Automatic Control, IEEE Transactions on},
  vol.~43, no.~8, pp. 1136--1142, 1998.

\bibitem{terra2014optimal}
M.~H. Terra, J.~P. Cerri, and J.~Y. Ishihara, ``Optimal robust linear quadratic
  regulator for systems subject to uncertainties,'' \emph{Automatic Control,
  IEEE Transactions on}, vol.~59, no.~9, pp. 2586--2591, 2014.

\bibitem{xie1993control}
L.~Xie and Y.~C. Soh, ``Control of uncertain discrete-time systems with
  guaranteed cost,'' in \emph{Decision and Control, 1993., Proceedings of the
  32nd IEEE Conference on}.\hskip 1em plus 0.5em minus 0.4em\relax IEEE, 1993,
  pp. 56--61.

\end{thebibliography}

\end{document}